\documentclass[preprint,12pt]{article}
\usepackage{amssymb}
\usepackage{amsmath}
\usepackage{hyperref}

\begin{document}

\centerline {\Large{\bf Solutions of the (3 + 1) dimensional 
 }}

\centerline{}

\centerline{\Large{\bf Charney - Obukhov equation for the ocean,     }}

\centerline{}

\centerline{\Large{\bf Part I: case of  ``separation'' of variables}}

\centerline{}

\centerline{\bf {V.A. Goloveshkin, A.G. Kudryavtsev, N.N. Myagkov}}

\centerline{}

\centerline{Institute of Applied Mechanics,}

\centerline{Russian Academy of Sciences, Moscow 125040, Russia}

\begin{abstract}

Exact solutions describing Rossby waves and vortices in ocean propagating along the zonal direction at a constant velocity are considered for the (3 + 1) -dimensional nonlinear Charney - Obukhov equation. In the first part of our work, we give examples of solutions of the Charney-Obukhov equation that satisfy the nonlinear boundary conditions of the ocean, for a special case that corresponds to the "separation" of variables.

\end{abstract}








In recent papers \cite {KudryavtsevPhL-2022}, \cite {KudryavtsevPhF-2022}, \cite {KudryavtsevPhF-2023}, exact solutions of the Charney–Obukhov equation for the ocean were obtained in the form of synoptic-scale waves and eddies against the background of a zonal flow depending on the vertical coordinate $z$. The zonal flow in these solutions is itself a solution to the Charney–Obukhov equation, while waves and vortices are not. These solutions are of physical interest, since they demonstrate the inseparable connection between the emerging motion of the ocean and the zonal flow.
It should be noted that for solutions describing drift in the zonal direction with a constant velocity $V$, Charney–Obukhov equation can be reduced to solving a linear Helmholtz-type equation. This reduction is well known in hydrodynamics and was used to obtain solutions in \cite {KudryavtsevPhL-2022}, \cite {KudryavtsevPhF-2022}, \cite {KudryavtsevPhF-2023}. To describe the processes in the ocean, the Charney–Obukhov equations are supplemented with nonlinear boundary conditions at the bottom and on the surface of the ocean. Taking into account the connection of the Charney-Obukhov equation with a linear equation, the Charney–Obukhov equation itself without taking into account the boundary conditions has a large set of solutions and a superposition of solutions. An interesting property of the considered model is the possibility of fulfilling nonlinear boundary conditions for a wide class of solutions of the Charney-Obukhov equation. In the first part of our work, we give examples of solutions of the Charney-Obukhov equation that satisfy the nonlinear boundary conditions of the ocean, for a special case that corresponds to the "separation" of variables.

It is well known that a large number of large scale circulations of the Earth's atmosphere and ocean is described by nonlinear quasi-geostrophic potential vorticity equation \eqref{eq1}, which is also known as the Charney - Obukhov equation. 
We consider the (3+1)-dimensional nonlinear Charney–Obukhov equation in the form
\begin{multline} \label{eq1}
{\frac {\partial ^{3}}{\partial {x}^{2}\partial t}}p \left( t,x,y,z
 \right) +{\frac {\partial ^{3}}{\partial {y}^{2}\partial t}}p \left( 
t,x,y,z \right) +{\frac {\partial ^{3}}{\partial {z}^{2}\partial t}}p
 \left( t,x,y,z \right) 
 \\
 + \left( {\frac {\partial }{\partial x}}p
 \left( t,x,y,z \right)  \right)  \left( {\frac {\partial ^{3}}{
\partial y\partial {x}^{2}}}p \left( t,x,y,z \right) +{\frac {
\partial ^{3}}{\partial {y}^{3}}}p \left( t,x,y,z \right) +{\frac {
\partial ^{3}}{\partial {z}^{2}\partial y}}p \left( t,x,y,z \right) 
 \right)
 \\
  - \left( {\frac {\partial }{\partial y}}p \left( t,x,y,z
 \right)  \right)  \left( {\frac {\partial ^{3}}{\partial {x}^{3}}}p
 \left( t,x,y,z \right) +{\frac {\partial ^{3}}{\partial {y}^{2}
\partial x}}p \left( t,x,y,z \right) +{\frac {\partial ^{3}}{\partial 
{z}^{2}\partial x}}p \left( t,x,y,z \right)  \right)
\\
 +\beta\,{\frac {
\partial }{\partial x}}p \left( t,x,y,z \right) =0
\end{multline}
Where $p$ is the dimensionless geostrophic stream function, bearing the sense of relative pressure perturbation; $\beta$  is the dimensionless meridional (northern) gradient of Coriolis parameter. As usual, we assume that the x-coordinate is east, the y-coordinate is north, and the z-coordinate is up.

We consider for equation \eqref{eq1} the boundary conditions for an ocean with a flat bottom and a rigid lid. Such boundary conditions are often used to model waves and vortices in the ocean, a detailed discussion of these boundary conditions is given, for example, in the book \cite {Pedlosky}. The boundary condition operator has the form
\begin{multline} \label{eq2}
{\frac {\partial ^{2}}{\partial z\partial t}}p \left( t,x,y,z \right) 
- \left( {\frac {\partial }{\partial y}}p \left( t,x,y,z \right) 
 \right) {\frac {\partial ^{2}}{\partial z\partial x}}p \left( t,x,y,z
 \right)
 \\
  + \left( {\frac {\partial }{\partial x}}p \left( t,x,y,z
 \right)  \right) {\frac {\partial ^{2}}{\partial z\partial y}}p
 \left( t,x,y,z \right) 
\end{multline}

We require that the boundary condition operator be equal to zero at $z=0$ and $z=H$, where  H is the ocean depth ( $0 \le z \le H$).

Solution 1.

The solution that satisfies the boundary conditions is
\begin{equation} \label{eq3}
F \left( x-Vt,y \right)
+ \left( \cos \left( {\frac {\pi \,nz}{H}} \right) P-{\frac {\beta\,{H
}^{2}}{{\pi }^{2}{n}^{2}}}-V \right) y
\end{equation}
where $P$ - arbitrary constant,  $n=1,2,3,...$, the function $F$ satisfies the equation
\begin{equation} \label{eq4}
{\frac {\partial ^{2}}{\partial {s}^{2}}}F \left( s,y \right) +{\frac 
{\partial ^{2}}{\partial {y}^{2}}}F \left( s,y \right) +\frac {{\pi }^{2}{n}^{2}}{{H}^{2}}F \left( s,y \right) = 0
\end{equation}

The boundary condition operator on this solution has the form
\begin{equation} \label{eq5}
- \left( {\frac {\partial }{\partial x}}F \left( x-Vt,y \right) 
 \right) \sin \left( {\frac {\pi \,nz}{H}} \right) \pi \,nP{H}^{-1}
\end{equation}

Solution 2.

The solution that satisfies the boundary conditions is
\begin{equation} \label{eq6}
\sin \left( {\it k_z}\,z \right) F \left( x-Vt,y \right) + \left( {
\frac {\beta\, \left( \sin \left( Kz \right) M+\cos \left( Kz \right) 
-1 \right) }{{K}^{2}}}-V \right) y
\end{equation}
where
${K}^{2}={{\it k_z}}^{2}+{{\it K_r}}^{2}$, ${\it k_z}, {\it K_r}$ - arbitrary constants, 
\begin{equation} \label{eq7}
M={\frac {{\it k_z}\,\cos \left( {\it k_z}\,H \right) \cos \left( KH
 \right) -{\it k_z}\,\cos \left( {\it k_z}\,H \right) +K\sin \left( {
\it k_z}\,H \right) \sin \left( KH \right) }{ K\sin \left( {\it k_z
}\,H \right) \cos \left( KH \right) -{\it k_z}\,\cos \left( {\it k_z}
\,H \right) \sin \left( KH \right)}}
\end{equation}
The function $F$ satisfies the equation
\begin{equation} \label{eq8}
{\frac {\partial ^{2}}{\partial {s}^{2}}}F \left( s,y \right) +{\frac 
{\partial ^{2}}{\partial {y}^{2}}}F \left( s,y \right) +{{\it K_r}}^{2}F \left( s,y \right) = 0
\end{equation}

The boundary condition operator on this solution has the form
\begin{multline} \label{eq9}
 \left( {\frac {\partial }{\partial x}}F \left( x-Vt,y
 \right)  \right) \beta {K}^{-2} \, \left( M \left( K\sin \left( {\it k_z}\,z
 \right) \cos \left( Kz \right) -{\it k_z}\cos \left( {\it k_z}\,z
 \right) \sin \left( Kz \right)  \right)
 \right.
\\ 
\left.
  -{\it k_z}\cos \left( {\it k_z}\,z
 \right) \cos \left( Kz \right) + {\it k_z} \cos \left( {\it k_z}\,z
 \right)-K\sin \left( {\it k_z}\,z \right) \sin \left( Kz
 \right)  \right) 
\end{multline}

Solution 3.

The solution that satisfies the boundary conditions is
\begin{equation} \label{eq10}
\cos \left( {\it k_z}\,z \right) F \left( x-Vt,y \right) + \left( {
\frac {\beta\, \left( \cos \left( Kz \right) M-1 \right) }{{K}^{2}}}-V
 \right) y
\end{equation}
where
${K}^{2}={{\it k_z}}^{2}+{{\it K_r}}^{2}$, ${\it k_z}, {\it K_r}$ - arbitrary constants, 
\begin{equation} \label{eq11}
M={\frac {{\it k_z}\,\sin \left( {\it k_z}\,H \right) }{{\it k_z}\sin \left( {
\it k_z}\,H \right) \cos \left( KH \right) -K\cos \left( 
{\it k_z}\,H \right) \sin \left( KH \right) }}
\end{equation}
The function $F$ satisfies the equation \eqref{eq8}.

The boundary condition operator on this solution has the form
\begin{multline} \label{eq12}
\left( {\frac {\partial }{\partial x}}F \left( x-Vt,y \right)  \right) \beta {K}^{-2}
\left( M \left( {\it k_z}\,\sin \left( {\it k_z}\,z \right) \cos \left( Kz \right) 
-K\cos \left( {\it k_z}\,z \right) \sin \left( Kz \right)  \right)
 \right.
\\ 
\left.
  - {\it k_z} \sin \left( {\it k_z}\,z \right)\right) 
\end{multline}

Solution 4.

The solution that satisfies the boundary conditions is
\begin{equation} \label{eq13}
\sinh \left( {\it k_z}\,z \right) F \left( x-Vt,y \right) + \left( {
\frac {\beta\, \left( \sin \left( Kz \right) M+\cos \left( Kz \right) 
-1 \right) }{{K}^{2}}}-V \right) y
\end{equation}
where
${K}^{2}=-{{\it k_z}}^{2}+{{\it K_r}}^{2}$, ${\it k_z}, {\it K_r}$ - arbitrary constants, ${\it k_z} < {\it K_r}$,  
\begin{equation} \label{eq14}
M={\frac {\cosh \left( {\it k_z}\,H \right) {\it k_z}\,\cos \left( KH
 \right) -\cosh \left( {\it k_z}\,H \right) {\it k_z}+\sinh \left( {
\it k_z}\,H \right) K\sin \left( KH \right) }{-\cosh \left( {\it k_z
}\,H \right) {\it k_z}\,\sin \left( KH \right) +K\sinh \left( {\it 
k_z}\,H \right) \cos \left( KH \right) }}
\end{equation}
The function $F$ satisfies the equation \eqref{eq8}.

The boundary condition operator on this solution has the form
\begin{multline} \label{eq15}
 \left( {\frac {\partial }{\partial x}}F \left( x-Vt,y
 \right)  \right) \beta  {K}^{-2} \,  \left( -\cosh \left( {\it k_z}\,z \right) {
\it k_z}\,\sin \left( Kz \right) M 
 \right.
\\ 
\left.
+\cosh \left( {\it k_z}\,z \right) {\it k_z}
-\cosh \left( {\it k_z}\,z
 \right) {\it k_z}\,\cos \left( Kz \right) 
  \right.
\\ 
\left.
 +\sinh \left( {\it k_z}\,z \right) K\cos \left( K
z \right) M-\sinh \left( {\it k_z}\,z \right) K\sin \left( Kz
 \right)  \right) 
\end{multline}

Solution 5.

The solution that satisfies the boundary conditions is
\begin{equation} \label{eq16}
\sinh \left( {\it k_z}\,z \right) F \left( x-Vt,y \right) + \left( {
\frac {\beta\, \left( \sinh \left( Kz \right) M-\cosh \left( Kz
 \right) +1 \right) }{{K}^{2}}}-V \right) y
\end{equation}
where
${K}^{2}={{\it k_z}}^{2}-{{\it K_r}}^{2}$, ${\it k_z}, {\it K_r}$ - arbitrary constants, ${\it k_z} > {\it K_r}$,  
\begin{equation} \label{eq17}
M={\frac {-\cosh \left( {\it k_z}\,H \right) {\it k_z}\,\cosh \left( K
H \right) +\cosh \left( {\it k_z}\,H \right) {\it k_z}+\sinh \left( 
{\it k_z}\,H \right) K\sinh \left( KH \right) }{-\cosh \left( {\it 
k_z}\,H \right) {\it k_z}\,\sinh \left( KH \right) +K\sinh \left( {
\it k_z}\,H \right) \cosh \left( KH \right) }}
\end{equation}
The function $F$ satisfies the equation \eqref{eq8}.

The boundary condition operator on this solution has the form
\begin{multline} \label{eq18}
  \left( {\frac {\partial }{\partial x}}F \left( x-Vt,y
 \right)  \right) \beta  {K}^{-2} \, \left( -\cosh \left( {\it k_z}\,z \right) {
\it k_z}\,\sinh \left( Kz \right) M
 \right.
\\ 
\left.
 -\cosh \left( {\it k_z}\,z \right) {\it k_z}
+\cosh \left( {\it k_z}\,z \right) {\it k_z}\,\cosh \left( Kz \right)
 \right.
\\ 
\left.
+\sinh \left( {\it k_z}\,z \right) K\cosh
 \left( Kz \right) M-\sinh \left( {\it k_z}\,z \right) K\sinh \left( 
Kz \right)  \right) 
\end{multline}

Solution 6.

The solution that satisfies the boundary conditions is
\begin{equation} \label{eq19}
\sinh \left( {\it k_z}\,z \right) F \left( x-Vt,y \right) + \left( -1
/2\,\beta\,{z}^{2}+1/2\,\beta\,HMz-V \right) y
\end{equation}
where
 ${\it k_z}$ - arbitrary constant,  
\begin{equation} \label{eq20}
M={\frac {2\,\sinh \left( {\it k_z}\,H \right) -\cosh \left( {\it k_z}
\,H \right) {\it k_z}\,H}{\sinh \left( {\it k_z}\,H \right) -\cosh
 \left( {\it k_z}\,H \right) {\it k_z}\,H}}
\end{equation}
The function $F$ satisfies the equation
\begin{equation} \label{eq21}
{\frac {\partial ^{2}}{\partial {s}^{2}}}F \left( s,y \right) +{\frac 
{\partial ^{2}}{\partial {y}^{2}}}F \left( s,y \right) +{{\it k_z}}^{2}F \left( s,y \right) = 0
\end{equation}

The boundary condition operator on this solution has the form
\begin{multline} \label{eq22}
-1/2\, \left( {\frac {\partial }{\partial x}}F \left( x-Vt,y \right) 
 \right) \beta\, \left( -\cosh \left( {\it k_z}\,z \right) {\it k_z}
\,{z}^{2}+\cosh \left( {\it k_z}\,z \right) {\it k_z}\,HMz
 \right.
\\ 
\left.
+2\,\sinh
 \left( {\it k_z}\,z \right) z-\sinh \left( {\it k_z}\,z \right) HM
 \right) 
\end{multline}

Solution 7.

The solution that satisfies the boundary conditions is
\begin{equation} \label{eq23}
\cosh \left( {\it k_z}\,z \right) F \left( x-Vt,y \right) + \left( {
\frac {\beta\, \left( \cos \left( Kz \right) M-1 \right) }{{K}^{2}}}-V
 \right) y
\end{equation}
where
${K}^{2}=-{{\it k_z}}^{2}+{{\it K_r}}^{2}$, ${\it k_z}, {\it K_r}$ - arbitrary constants, ${\it k_z} < {\it K_r}$,  
\begin{equation} \label{eq24}
M={\frac {\sinh \left( {\it k_z}\,H \right) {\it k_z}}{\sinh \left( {
\it k_z}\,H \right) {\it k_z}\,\cos \left( KH \right) +K\cosh
 \left( {\it k_z}\,H \right) \sin \left( KH \right) }}
\end{equation}
The function $F$ satisfies the equation \eqref{eq8}.

The boundary condition operator on this solution has the form
\begin{multline} \label{eq25}
 - \left( {\frac {\partial }{\partial x}}F \left( x-Vt,y
 \right)  \right) \beta  {K}^{-2} \, \left( \sinh \left( {\it k_z}\,z \right) {
\it k_z}\,\cos \left( Kz \right) M-\sinh \left( {\it k_z}\,z
 \right) {\it k_z}
  \right.
\\ 
\left.
 +\cosh \left( {\it k_z}\,z \right) K\sin \left( Kz
 \right) M \right) 
\end{multline}

Solution 8.

The solution that satisfies the boundary conditions is
\begin{equation} \label{eq26}
\cosh \left( {\it k_z}\,z \right) F \left( x-Vt,y \right) + \left( {
\frac {\beta\, \left( \cosh \left( Kz \right) M+1 \right) }{{K}^{2}}}-
V \right) y
\end{equation}
where
${K}^{2}={{\it k_z}}^{2}-{{\it K_r}}^{2}$, ${\it k_z}, {\it K_r}$ - arbitrary constants, ${\it k_z} > {\it K_r}$,  
\begin{equation} \label{eq27}
M={\frac {\sinh \left( {\it k_z}\,H \right) {\it k_z}}{-\sinh \left( {
\it k_z}\,H \right) {\it k_z}\,\cosh \left( KH \right) +K\cosh
 \left( {\it k_z}\,H \right) \sinh \left( KH \right) }}
\end{equation}
The function $F$ satisfies the equation \eqref{eq8}.

The boundary condition operator on this solution has the form
\begin{multline} \label{eq28}
\left( {\frac {\partial }{\partial x}}F \left( x-Vt,y
 \right)  \right) \beta {K}^{-2} \, \left( -\sinh \left( {\it k_z}\,z \right) {
\it k_z}\,\cosh \left( Kz \right) M-\sinh \left( {\it k_z}\,z
 \right) {\it k_z}
  \right.
\\ 
\left.
 +\cosh \left( {\it k_z}\,z \right) K\sinh \left( K
z \right) M \right) 
\end{multline}

Solution 9.

The solution that satisfies the boundary conditions is
\begin{equation} \label{eq29}
\cosh \left( {\it k_z}\,z \right) F \left( x-Vt,y \right) + \left( -1
/2\,\beta\,{z}^{2}+1/2\,\beta\,{H}^{2}M-V \right) y
\end{equation}
where
 ${\it k_z}$ - arbitrary constant,  
\begin{equation} \label{eq30}
M=1-2\,{\frac {\cosh \left( {\it k_z}\,H \right) }{\sinh \left( {\it 
k_z}\,H \right) H{\it k_z}}}
\end{equation}
The function $F$ satisfies the equation \eqref{eq21}.

The boundary condition operator on this solution has the form
\begin{multline} \label{eq31}
-1/2\, \left( {\frac {\partial }{\partial x}}F \left( x-Vt,y \right) 
 \right) \beta\, \left( -\sinh \left( {\it k_z}\,z \right) {\it k_z}
\,{z}^{2}+\sinh \left( {\it k_z}\,z \right) {\it k_z}\,{H}^{2}M
 \right.
\\ 
\left.
+2\,
\cosh \left( {\it k_z}\,z \right) z \right) 
\end{multline}

All presented solutions have the form 
\begin{equation} \label{eq32}
f(z)F \left( x-Vt,y \right) +U(z;V)y
\end{equation}
where $U(z;V)$ is the velocity of the steady zonal background flow which depends on the velocity $V$. Thus, the first term of the formula \eqref{eq32} can be represented as the product of the horizontal coordinate function and the vertical coordinate function, this is a case of "separation" of variables.

\end{document}